\documentclass[12pt]{article}
\usepackage{epic,eepic,epsf,epsfig}
\usepackage{amsmath}
\usepackage{amssymb}
\usepackage{amsbsy}

\usepackage{amsfonts}

\begin{document}

\baselineskip 16pt

\title{ \vspace{-1.2cm}
On the widths of finite groups (I)}

\author{Wujie Shi\\\
\\
{\small{ Department of Mathematics, Chongqing University
of Arts and Sciences,}}\\
{\small{Chongqing 402160, P.R.China}}\\
{\small{ and}}\\
{\small{ School of Mathematics, Suzhou University, Suzhou
215006, P.R.China}}\\
{\small{ Email: shiwujie@outlook.com}}\\}

\date{ }

\maketitle \vspace{-.8cm}

\textbf{Abstract:} Let $G$ be a finite group, $\pi(G)$ be the set of prime divisors dividing the order of G and $\pi_{e}(G)$ (spectrum) denote the set of element orders of $G$.  We define $w_{o}(G)=
|\pi(G)|$ the width of order of G and $w_{s}(G)= max \{|\pi(k)| \mid  k \in \pi_{e}(G)\}$ the width of spectrum of $G$. In this paper, we discuss the cases of $w_{o}(G)$ and $w_{s}(G)$ are small,  prove several new results and give a survey about the two widths of groups. This article is a part of revised based on the author's plenary talk at ``2020 Ural Workshop on Group Theory and Combinatorics".



\

 Let $G$ be a finite group,  $\pi(G)$ be the set of prime divisors of $|G|$, the order of $G$. Let $w_{o}(G)=
|\pi(G)|$ denote the width of order of $G$, that is, the number of different prime divisors of $|G|$. Let $\pi_{e}(G)$ denote the set of element orders of $G$, it is called the spectrum of $G$. Denote $w_{s}(G)= max \{|\pi(k)| \mid  k \in \pi_{e}(G)\}$ the width of spectrum of $G$. In this paper, we discuss the cases of $w_{o}(G)$ and $w_{s}(G)$ are small,  give a survey about the two widths of groups, prove several new results and pose some problems.

\textbf{Finite groups with smaller widths of order}

1. A finite group with $w_{o}(G)=1$ is a finite $p$ group. Although the group order has only one prime factor, the structure is still very complicated. For the recent review literature, please refer to Zhang Qinhai and An Lijian's monographs[1, 2]. Looking at a special ``simple" situation, $\pi_{e}(G) = \{1, 3\}$ (``prime element" 3 groups), which the structure is still complicated, and so far, there is no clear classification.

2. The finite group with $w_{o}(G)=2$, namely the group of order $p^{a}q^{b}$, is a special kind of solvable group by the famous $p^{a}q^{b}$ theorem. W. Burnside gave the first proof in 1904 using representation theory [3], H. Bender gave the abstract group proof of this theorem in 1972 [4].

3. For the case of $w_{o}(G)=3$, they can be a nonsolvable groups, and their chief series can contain a nonabelian simple group with $w_{o}(G)=3$. Ref. [5] gives a simple group with $w_{o}(G)=3$. They are isomorphic to one of the following simple groups: $A_{5}, L_{2}(7), L_{2}(8), A_{6}, L_{2}(17), L_{3}(3), U_{3}(3)$ or $U_{4}(2)$. D. Gorenstein called the above 8 simple groups a simple $K_{3}$ group ($w_{o}(G)=3$), and pointed out ``that classification is obtained only as a corollary of the complete classification of all finite simple groups" (CFSG)[6].

4. Following the simple $K_{3}$ group, we consider the problem of the simple $K_{4}$ group using CFSG. The refs. [$7^{[Note]}$, 8, 9] successively gives all the simple $K_{4}$ groups, which are one of the following groups: $A_{n}, n = 7, 8, 9, 10$; $M_{11}, M_{12}, J_{2}$; $L_{2}(q), q = 16, 25, 49, 81$; $L_{3}(q), q = 4, 5, 7, 8, 17$; $L_{4}(3)$; $O_{5}(q), q = 4, 5, 7, 9$; $O_{7}(2), O_{8}^{+}(2), G_{2}(3)$; $U_{3}(q), q = 4, 5, 7, 8, 9$; $U_{4}(3)$; $U_{5}(2); ^{3}D_{4}(2)$; $^{2}F_{4}(2)'$ ; $Sz(8), Sz(32)$; and
$L_{2}(r)$, $r$ being prime and satisfying the following equation:

\begin{center}                                          $r^{2}-1=2^{a}3^{b}u^{c}$    \  \   \  \  \  \ \  \   \  \   (1)
\end{center}

where $a\geq1$, $b \geq1$, $c\geq1$, $u$ prime, $u >3$;

$L_{2}(2^{m})$ and satisfying the following equations:

   \begin{center}                                         $2^{m}-1=u,2^{m}+1=3t^{b}$            \  \   \  \  \  \ \  \   \  \  (2)
    \end{center}

where $m\geq1$, $u, t$ primes, $t >3, b\geq1$.

$L_{2}(3^{m})$ and satisfying the following equations:

\begin{center}                                           $3^{m}+1=4t,3^{m}-1=2u^{c}$     \  \   \  \  \  \ \  \   \  \    (3)

\end{center}
or

\begin{center}                                           $3^{m}+1=4t^{b},3^{m}-1=2u$     \  \   \  \  \  \ \  \   \  \       (4)
\end{center}

where $m \geq1$, $u, t$ odd primes, $c \geq1, b \geq1$.

After the publication of ref. [8], [10] discussed the above equations (1)-(4), and gave all their solutions using algebraic number theory. However, [11] pointed out: the proof of [10 ] cited article is not correct, so the conclusion of [10] is incorrect. Using the computer program ref. [12] got out, if the largest prime divisor of the group order is less than 1060, then there are only 101 simple $K_{4}$ groups. Ref. [13] proved that, in the above equations (2), (4), if $b > 1$, then the equations have no solution, while in the equations (3), when $c > 1$, there are only solutions $(u, t, m, c) = (11, 61, 5, 2)$. [14] proved that for equation (1), if $c > 1$, there are only solutions $(r, u, a, b, c) = (97, 7, 6, 1, 2)$ or $(r, u, a, b, c) = (577, 17, 7, 2, 2)$. Furthermore, Ref. [14] believes that the proof of the infinity of the number of simple $K_{4}$ groups is a difficult problem, which is a special case of the unresolved Dickson conjecture [15]. It can be attributed to whether there are infinite prime number pairs in the following number pairs: $(x, 3x-2)$, $(x, 2x + 1)$, $(x, 4x + 1)$, $(x, 6x-1)$, $(x, 6x + 1)$, $(x, 2x-1)$, $(x, 4x-1)$, $(x, 8x-1)$, where $x $ prime. If there are infinitely many prime numbers in one of the number pairs, it proves that the number of simple $K_{4}$ groups is infinite. On the contrary, if the second number in the above-mentioned number pair is only a finite prime numbers, can get to contradictions from this assumption?

\textbf{Conjecture 1.} There exist infinite simple $K_{4}$ groups.

5. Following ref. [8], [16] gave the simple groups with $w_{o}(G) = 5, 6$, namely the classification of simple $K_{5}$ groups and simple $K_{6}$ groups. Similarly, there are also such problem, the number of such groups is finite or infinite?

6. For the special number $|G|$ we can determine the group $G$, for example, $|G| = p$, $p$ is a prime. That is, there is only one group of order $p$ which is an abelian simple group. Let $G$ be a finite group. Denote $|G| = n$ and $f(n)$ the number of distinct groups $G$ with
$|G| = n$. How many different groups with the given order exist? The question was first posed by Cayley in 1854. The first study of the values of $f(n)$ seems to be in the work of $G$. A. Miller in the 1920's. O. H$\ddot{o}$lder determines all groups of square-free order and gives an explicit formula for their number. Many specialists of group theory re-searched also it, as J.G. Thompson, G. Higman, C.C. Sims and others. We may refer to [16] and the references listed in [17].

 $f(n) = 1$ if and only if $(n, \varphi(n)) = 1$, where $\varphi(n)$ is Euler function of $n$. Also the solutions for $f(n) = 2, 3$ and 4 is found [18].
 
Discuss the inverse function $f^{-1}$, for given $k$, which $n$ satisfies $f(n) = k$ ? That is, if $n$ is an integer such that the number of isomorphism classes of groups of order $n$ is $k$or a given $k$, then $n$ is said to be a solution of the equation $f(n) = k $. It is proved that if $1\leq k\leq 1000$, then there are infinitely many solutions for $f(n) = k$ ( see [18, 19]).

\textbf{Conjecture 2.}  For any $k\in Z^{+}$, there exist an integer $n$, such that $f(n) = k$.

7. If assume $G$ is simple, then the number of distinct simple groups $G$ with $|G| = n$,  $f_{simple}(n) = 0, 1$ or 2 [20, 21].

Let $G$ be a simple group and $|G| = |M|$, where $M$ is a known finite simple group. Then one of the following conclusions holds.

(1)	 If $|M| = |A_{8}| = |L_{3}(4)|$, then $G\simeq A_{8}$ or $L_{3}(4)$;

(2)	If $|M| = |B_{n}(q)| = |C_{n}(q)|$ where $n \geq 3$ and $q$ is odd , then $G\simeq B_{n}(q)$ or $C_{n}(q)$;

(3)	If $|M|$ is not the above case (1) or (2), then $G\simeq M$.

8. By inspecting the orders of all the known finite simple groups we have the following statement [21]:

Let $G$ be a finite non-abelian simple group. Then $G$ contains a proper subgroup $H$ such that $|G| < |H|^{2}.$

Let $G$ be a finite group and $|G|\neq p$, $p $ prime.Then $G$ contains a proper subgroup $H$ such that $|G|\leq |H|^{2}$.

9. Before we discuss the widths of spectrum $w_{s}(G)$, we comparing the two numbers $w_{o}(G)$ and $w_{s}(G)$ firstly.

\textbf{Theorem 1.} Let $G$ be a finite non-abelian simple group. Then  $w_{s}(G) < w_{o}(G)$.

\textbf{Proof.}  Since $w_{s}(G) = max \{ |\pi(k)| | k \in \pi_{e}(G) \}$, the set $\pi_{e}(G)$ defines the prime graph (Gruenberg-Kegel graph) $\Gamma(G)$ of $G$; in this graph the vertex set is $\pi(G)$, and distinct vertices $p$ and $q$ are adjacent if and only if $pq \in \pi_{e}(G)$. The concept of prime graph appeared in the unpublished manuscript [22] by K. Gruenberg and O. Kegel. This result was published later ref. [23]. [23] obtained this description for all simple groups except simple groups of Lie type in characteristic 2. The case of simple groups of Lie type in characteristic 2 were described by [24], later this result was obtained independently by [25, 26]. Unfortunately, all the papers contain some inaccuracies. Most of these inaccuracies was corrected in [27], and then the corrections were finished in [28]. Let G be a finite non-abelian simple group and denote the connected components of the prime graph by $\{\pi_{i}, i = 1, 2, \ldots, t \}$, $t$ the number of the connected component of $G$, and if the order of $G$ is even, denote the component containing 2 by $\pi_{1}.$

    Obviously, $w_{s}(G)\leq w_{o}(G)$ by the definitions. If $t > 1$, the conclusion is hold. So, we consider only the case in which the prime graph of $G$ is connected. That is, the number of component is one. From [23] and [24], we consider the following simple groups:

(1)	$A_{n},  n \geq 5, n \neq p, p +1,  p +2, p$ prime;

(2)	$E_{7}(q), q\neq 3$;

(3)  classical groups $A_{l},  B_{l},  C_{l},  D_{l},  ^{2}A_{l},  ^{2}D_{l}$ satisfy some conditions.

Now we prove this theorem according the above three cases:

(1)	For the case of alternating groups $A_{m}$, we have the following lemma:

\textbf{Lemma 1.} ([29, Lemma 4]) In $A_{m}$ , there is an element of order $n = p_{1}^{a}p_{2}^{b}\cdots p_{s}^{c}$, where  $p_{1},p_{2},\ldots, p_{s}$ are distinct primes and $a, b, \ldots , c$ are naturals, if and only if  $ p_{1}^{a}+p_{2}^{b}+\cdots +p_{s}^{c}\leq m$ for odd   $n$ or  $ p_{1}^{a}+p_{2}^{b}+\cdots +p_{s}^{c}\leq m-2$ for even $n$.

\textbf{Proof of Case (1).} If $w_{s}(G)  = w_{o}(G)$, then G contains the elements of order $p_{1}p_{2}\cdots p_{s}$ , where $w_{o}(G) =  |\pi(G)| = s$. Next we prove that  $ p_{1}+p_{2}+\cdots +p_{s} > m$, where $|A_{m}| = m!/2$. For the case of $m < 12$ we may check it immediately. For example $m = 10$ and $|A_{10}| = (1/2) (10 !)$, and $|\pi(A_{10})| = 5, 2+3+5+7 > 10$.

For $m \geq 12$, there is two primes $n_{1}, n_{2}$  between $[m/2]$ and $m$ by the prime number theorem (PNT). Again, $n_{1} + n_{2} > m$, we get contradiction. In fact, ref. [30] showed there is a prime in the interval $[3k, 4k]$ for all $k \geq 2$. Thus, there are two primes in the interval $[2k, 4k]$ for all $k \geq 3$.

(2) The case of $E_{7}(q), q\neq 3$; we list the following lemma:

\textbf{Lemma 2.} ([31]) Let $k \geq 2$ an integer and $p$ a prime.

(a) There exists a primitive prime divisor of $p^{k}- 1$ unless $k = 2$ and $p$ is a Mersenne prime or $(p, k) = (2, 6).$

(b) If $r$ is a primitive prime divisor of $p^{k} - 1$, then $r-1\equiv 0 (mod \ k)$. In particular, $r\geq k + 1.$

\textbf{Lemma 3.} ([31]) Let $G \simeq E_{7} (q), q = p^{n}$. Then

(a) $|G| = (1/d) q^{63}(q^{18} -1)(q^{14} -1)(q^{12} -1)(q^{10} -1)(q^{8} -1)(q^{6} -1)(q^{2} -1), d = (2, q -1).$

(b) Let $a, b \in E_{7} (q), q = p^{n}$. $|a| = p$, and $ab\neq ba.$ Then $T = T1\cup T2\subseteq \{ |b| | ab \neq ba\}$ and $T1\cap T2 = \varnothing$, where $T = \pi ^{*'} ((q^{6}-q^{3}+1)(q^{7}+1)/(q+1))$, $T1 = \pi^{*}(q^{6}-q^{3}+1)$,  $T2 =
     \pi^{'}((q^{7}+1)/(q+1))$, and $\pi^{*'}(k)$ denotes the set of prime factors different from 3 and    7 in $k$, $\pi^{*}(k)$ denotes the set of prime factors different from 3 in $k$, and $\pi^{'}(k)$ denotes the set of prime factors different from 7 in $k$.

\textbf{Proof of Case (2).} If $w_{s}(G) = w_{o}(G)$, then $G$ contains the elements of order $p_{1}p_{2}\cdots p_{s}$ , where $w_{o}(G) =  |\pi(G)| = s$, then $G$ contains the elements of order $pp_{i}$ . It is contrary to that $T\neq \varnothing$(in Lemma 3 (b) and Lemma 2).

(3)  The case of classical groups $A_{l},  B_{l},  C_{l},  D_{l},  ^{2}A_{l},  ^{2}D_{l}$ satisfy some conditions.
In this case we consider the orders of the centralizer of involution in $G$. If $w_{s}(G)  = w_{o}(G)$, then $G$ contains the elements of order $p_{1}p_{2}\cdots p_{s}$ , where $w_{o}(G) =  |\pi(G)| = s$. That is, $G$ contains the elements of order 2, and the order of its centralizer covers all prime factors of$G|$. Next, we prove that it is impossible from ref. [23].

\textbf{Proof of Case (3).} The case of $G\simeq A_{l}(q), q = p^{n}.$ Then

\begin{center} $|G| = |A_{l}(q)| = (1/d) q^{l(l+1)/2} (q^{2}-1) (q^{3}-1)\cdots(q^{l+1}-1)$,  $d = (l+1, q-1)$,
\end{center}

and there are centralizers of type $A_{r}\times A_{s} , r + s = l - 1$ except $l = 1$. All the factors dividing these orders of centralizers are contained in the centralizer of type $A_{0} \times A_{l-1}$. Thus there is a prime $p_{j} , p_{j} | (q^{l+1}-1)$ and $p_{j}\nmid q^{l(l+1)/2} (q^{2}-1) (q^{3}-1)\cdots (q^{l}-1)$ except $l+1 = 2$ and $q$ is a Mersenne prime or $(q, l+1) = (2, 6)$ by Lemma 2. If $q = 2, l = 5$, we have $|A_{5}(2)| = 215\cdot34\cdot5\cdot72\cdot31$, and $\pi_{e}(G) = \pi_{e}(A_{5}(2))$ =$\{1, 2, \cdots, 10, 12, 14, 15, 21, 28, 30, 31, 63\}$. In this case, $w_{s}(G) = 3 < w_{o}(G) = 5$. If $l = 1$, then the group $A_{1}(q)$ is isomorphic to $PSL(2, q)$ and the number of prime graph components more than $1$, thus $w_{s}(G) < w_{o}(G)$.  For the cases of $B_{l},  C_{l},  D_{l},  ^{2}A_{l},  ^{2}D_{l}$ , we can check the conclusion hold using the similar verification.

Since $w_{s}(G) < w_{o}(G)$ for all finite simple groups we consider the difference of two widths $d = w_{o}(G) - w_{s}(G), d \geq 1 $. It is easy to verify that no alternating groups and no sporadic simple groups whose difference of two widths equal to 1.

The following simple groups of Lie type whose difference of two widths $d =1$:

$L_{3}(3), U_{3}(3), U_{4}(2)$ (for all $w_{o}(G) = 3$);

$S_{4}(5), S_{4}(7), S_{4}(9), U_{3}(9)$ (for all $w_{o}(G) = 4$).

\textbf{Question 1.} Classify all finite simple groups whose difference of two widths equal to one.


\begin{thebibliography}{9}
\begin{small}

\bibitem{ZA} Zhang Qinhai, An Lijian, The Structure of Finite $p$-Groups (Vol. 1), Beijing: Sci-ence Press, 2017. (Chinese)
\bibitem{ZAL}  Zhang Qinhai, An Lijian, The Structure of Finite $p$-Groups (Vol. 2), Beijing: Sci-ence Press, 2017. (Chinese)
\bibitem{WB} W. Burnside, On groups of order $p^{\alpha}q^{\beta}$, London Math. Soc. s2-1 (1904), no.1, 388-392.
\bibitem{HB} H. Bender, A group theoretic proof of Burnside's $p^{a}q^{b}$-theorem, Math. Z. 126 (1972), 327-338.
\bibitem{MH} M. Herzog, On finite simple groups of order divisible by three primes only, J. Algebra, 10 (1968), 383-388.
\bibitem{DG}  D. Gorenstein, Finite Simple Groups, Plenum Press, New York and London, 1982.
\bibitem{WX}  Wang Xiaofeng, Classification theorem of $K_{4}$ simple groups, Chinese Science Bul-letin, no.14 (1990), 1117-1118. (Chinese)

         [Note] In [7] there are the conclusions without specific proofs, some omissions (such as $L_{3}(8), U_{3}(7)$), and repetitions (such as $O_{7}(2)$ and $S_{6}(2)$ are same group), in addition, Suzuki series simple groups can be completely determined.
\bibitem{SW}  Shi Wujie, On simple $K_{4}$-groups, Chinese Science Bulletin, no.17 (1991), 1281-1283. (Chinese)
\bibitem{BW}  B. Huppert, W. Lempken, Simple groups of order divisible by at most four primes, Proc. of the F. Scorina Gomel State University, No. 3(16) (2000), 64-75.
\bibitem{LX}  Le Maohua, Xu Guangshan, Several Diophantine equation problems of simple $K_{4}$-groups, Science in China (Ser. A), 26(1996), 769-773. (Chinese)
\bibitem{YP}  Yuan Pingzhi, Several unsolved Diophantine equations, J. Math. Res. Exposition 20 (2000), no.4, 627-628.  (Chinese)
\bibitem{DH}  Deng Huiwen, On the number of simple $K_{4}$-groups, Journal of Southwest Uni-versity (Natural Science Edition), 23:4 (1998), 375-378.(Chinese)
\bibitem{YZM}  Y. Bugeaud, Z. Cao, M. Mignotte, On simple $K_{4}$-groups,  J. Algebra 241 (2001), 658-668.
\bibitem{SZW}  Shaohua Zhang, Wujie Shi, On the Number of Simple $K_{4}$-Groups, Bulletin of the Iranian Mathematical Society,
\bibitem{LED}  L.E. Dickson, A new extension of Dirichlet's theorem on prime numbers. Messenger Math. 33(1904), 155-161.
\bibitem{JI}  A. Jafarzadeh, A. Iranmanesh, On simple $K_{n}$-groups for $n = 5, 6$, Groups St An-drews 2005, London Mathematical Society Lecture Notes Series (340), 517-526.
\bibitem{LP}  L. Pyber, Enumerating finite groups of given order, Ann. of Math.137(1993), 203-220.
\bibitem{GW}  Guimin. Wei, Enumeration Formulae for the number of finite groups and their application, Southeast Asian Bull. Math., 22 (1998), 93-102.
\bibitem{ZY}  Zhang Yuanbiao, Groups with less than 6000 isomorphism classes, Journal of Southwest University (Natural Science Edition), 29:9 (2007), 20-24.(Chinese)
\bibitem{WRRD}  W. Kimmerle, R. Lyons, R. Sandling and D.N. Teague, Composition factors from the groups ring and Artin's theorem on orders of simple groups, Proc. London Math. Soc. (3) 60 (1990), 89-122.
\bibitem{WJ}  Wujie Shi, On the orders of the finite simple groups, Chinese Science Bulletin, 38:4 (1993), 296-298. (Chinese)
\bibitem{GK}  K. W. Gruenberg, O. Kegel, Unpublished manuscript, 1975.
\bibitem{JS}  J. S. Williams, Prime graph components of finite groups, J. Algebra 69 (1981), 487-513.
\bibitem{AS}  A. S. Kondrat'ev, Prime graph components of finite simple groups, Math. USSR Sb. 67 (1990), 235-247.
\bibitem{IY}  N. Iiyori, H. Yamaki, Prime graph components of the simple groups of Lie type over the fields of even characteristic, Proc. Japan Acad. Ser. A. Math. Sci., 67:3 (1991), 82-83.
\bibitem{NIY}  N. Iiyori, H. Yamaki, Prime graph components of the simple groups of Lie type over the fields of even characteristic, J. Algebra, 155:2 (1993), 335-343; Corrigenda, J. Algebra, 181:2 (1996), 659.
\bibitem{KM}  A. S. Kondrat'ev and V. D. Mazurov, Recognizibility of alternating groups of prime degree by their element orders, Siberian Math. J. 41:2 (2000), 359-369.
\bibitem{ASK}  A. S. Kondrat'ev, Gruenberg-Kegel graph of a finite group and its applications, Algebra and Linear Optimization, Proc. of Internation in honor of 90th Birthday S. N. Chernikov, UB RAS: Ekaterinburg, 2002, 141-158 (in Russian).
\bibitem{ZM} A. V. Zavarnitsin and V. D. Mazurov, Element orders in coverings of symmetric and alternating groups, Algebra and Logic, 38(1999), 159-170.
\bibitem{AL} A. Loo, On the primes in the interval $[3n, 4n]$. (English summary) Int. J. Contemp. Math. Sci. 6 (2011), no. 37- 40, 1871-1882.
\bibitem{KZ} K. Zsigmondy, Zur Theorie der Potenzreste,  Monatsh. Math. Phys, 3 (1892), 265-284.
\bibitem{WJS} Wujie Shi, The pure quantitative characterization of finite simple groups (I), 4:3(1994), 316-326.


\end{small}

\end{thebibliography}
\end{document}